\pdfoutput=1
\documentclass[journal]{IEEEtran}
\usepackage{amsfonts}
\usepackage{bm}
\usepackage{amsmath,amsfonts,amsthm,amssymb}
\usepackage{setspace}
\usepackage{Tabbing}
\usepackage{fancyhdr}
\usepackage{lastpage}
\usepackage{extramarks}
\usepackage{chngpage}
\usepackage{algorithmic}
\usepackage{soul,color}
\usepackage{epsfig}
\usepackage{graphicx,float,wrapfig,subfigure}
\usepackage{dsfont}
\usepackage{longtable}
\usepackage{wrapfig}
\usepackage{stfloats}
\usepackage{cite}
\usepackage{graphicx}
\usepackage{array}
\usepackage{multirow}{\tiny }
\usepackage{multicol}
\usepackage{colortbl}
\usepackage{tabularx}
\usepackage{mdwmath}
\usepackage{mdwtab}
\usepackage{color}
\usepackage{verbatim}
\usepackage{amsmath}
\usepackage{tikz}
\usetikzlibrary{calc}
\usepackage{flowchart}
\usetikzlibrary{shapes.geometric, arrows}
\usepackage{overpic}
\usepackage{epstopdf}
\usepackage{url}
\usepackage[colorlinks,linkcolor=black,anchorcolor=black,citecolor=black,urlcolor=black]{hyperref} 

\usepackage{makecell}
\usepackage{textcomp}
\usepackage[numbers,sort&compress]{natbib}
\usepackage{algorithm}
\usepackage{algorithmic}

\newtheorem{remark}{Remark}

\usepackage{color}
\makeatletter %
\let\myorg@bibitem\bibitem
\def\bibitem#1#2\par{%
	\@ifundefined{bibitem@#1}{%
		\myorg@bibitem{#1}#2\par
	}{%
		\begingroup
		\color{\csname bibitem@#1\endcsname}%
		\myorg@bibitem{#1}#2\par
		\endgroup
	}%
}

\makeatother %

\allowdisplaybreaks

\usepackage{footnote}
\makesavenoteenv{tabular}
\makesavenoteenv{table}

\begin{document}

\title{
	\vspace{-6mm}
	Topology-free optimal power dispatch for distribution network considering security constraints and flexible building thermal inertia 
	}

\author{
Ge~Chen,~\IEEEmembership{Student Member,~IEEE,}
Hongcai~Zhang,~\IEEEmembership{Member,~IEEE,}
Ningyi~Dai,~\IEEEmembership{Senior Member,~IEEE,}
and~Yonghua~Song,~\IEEEmembership{Fellow,~IEEE}
\vspace{-8mm}

\thanks{	
	This work is supported in part by the Science and Technology Development Fund, Macau SAR (File no. 0137/2019/A3). G. Chen, H. Zhang, N. Dai, and Y. Song are with the State Key Laboratory of Internet of Things for Smart City and Department of Electrical and Computer Engineering, University of Macau, Macao, 999078 China (email: hczhang@um.edu.mo).
	}
}

\maketitle

\begin{abstract}
		With the increasing integration of distributed PV generation, the distribution network requires more and more flexibility to achieve the security-constrained optimal power dispatch. 
		However, the conventional flexibility sources usually require additional investment cost for equipment. Moreover, involving the security constraints is very challenging	due to the requirements of accurate network model that may be unavailable in practice. This paper addresses the aforementioned challenge by proposing a topology-free optimal power dispatch framework for distribution networks. It utilizes building thermal inertia to provide flexibility to avoid additional investment. To guarantee the operation safety, a multi-layer perception (MLP) is trained based on historical operational data and then reformulated as mixed-integer  constraints to replace the inexplicit original security constraints. 
		Numerical results confirm that the proposed framework can derive a feasible and optimal strategy without any topology information.

\end{abstract}
\begin{IEEEkeywords}
	Demand response, building thermal inertia, optimization, security constraints, neural networks.
\end{IEEEkeywords}

\vspace{-2mm}
\section{Introduction}
\IEEEPARstart{I}{n} recent years, the investment in PV generation grows rapidly in response to the challenges of energy crisis and $\text{CO}_2$ emission. In 2017, the cumulative installed PV generation capacity reaches to 130.25 GW in China \cite{tang2018solar}. Because that PV generation is intermittent, stochastic and often uncontrollable, its increasing capacity may deteriorate operational efficiency and security of power distribution networks, e.g. leading to voltage instability and line current violation. In order to enhance operational efficiency and security in distribution networks with high penetration of PV generation, more flexibility sources and smarter power scheduling are required \cite{mohandes2019review}. 

Many scholars focus on coordinating distribution networks with distributed energy storage to improve system efficiency and security. References \cite{guo2018data1,guo2018data2} combined the battery with grid to reduce the PV curtailment. The security constraints, such as the line current and bus voltage limitations, were introduced based on the power flow equations. In \cite{biswas2017optimal}, the grid was coordinated with a pumped hydro system to store the excess renewable generation. Reference \cite{kikusato2019electric} treated electric vehicles as distributed energy storage systems and shifted their
charging period to effectively utilize PV power. In \cite{zhang2019development}, an electricity-gas integrated energy system was proposed to enhance the utilization of renewable energy and reduce operational cost.

Recently, an increasing number of works utilized building thermal inertia, such as controlling heating, ventilation, and air conditioning, to shift the peak loads and provide flexibility for distribution network. Compared with energy storage systems, thermal inertia does not require any additional devices, so the investment cost can be reduced.
Reference \cite{li2020collaborative} utilized thermal inertia of water in pipelines to reduce the operation cost of integrated electricity and district heating systems. In \cite{chen2020optimal}, the thermal inertia of pipeline network was treated as thermal storage systems to realize the optimal power dispatch for district cooling systems. Reference \cite{xue2019coordinated} considered the flexibility of heat sources, pipelines, and loads to avoid unnecessary renewable power curtailment. Reference \cite{zhang2019building} developed an approach to quantify the demand response flexibility provided by the building thermal inertia. In references \cite{shi2019thermostatic,adhikari2019heuristic}, aggregation models for heating, ventilation and air conditioning systems were established for providing regulation service.

Although the above papers have made great progress on the power dispatch of distribution network with high-penetration renewable energy, there are still significant gaps to fill. References \cite{li2020collaborative,zhang2019building, chen2020optimal,shi2019thermostatic,adhikari2019heuristic} utilized the demand-side sources to avoid the additional investment, but the power flow security constraints are ignored. References \cite{guo2018data1,guo2018data2,biswas2017optimal,kikusato2019electric,zhang2019development} built the security constraints via the power flow equations, but they require accurate distribution network model that may be unknown in many practical cases. Recently, deep reinforcement learning was employed to optimize the power dispatch \cite{lu2019incentive}. However, unlike the white-box models which can be solved by off-the-shelf solvers efficiently, those learning-based models are time-consuming and can not guarantee global optimality.

To overcome the aforementioned challenges, we propose a topology-free optimal power dispatch framework for distribution networks. This method advances the published literature in the following aspects: 
\begin{enumerate}
    \item We propose a topology-free  optimal power dispatch method for distribution networks. It trains a multi-layer perceptron (MLP) neural network based on historical data to describe the security constraints of the distribution network with unknown topology and parameters.
	\item We utilize buildings' thermal inertia as distributed energy storage systems to provide flexibility. This can significantly enhance the power dispatch performance and reduce PV curtailment. 
	\item A white-box optimal model is proposed by equivalently reformulating the intractable black-box MLP as mixed-integer linear constraints so that the global optimal strategy can be found efficiently by off-the-shelf solvers.
\end{enumerate}
The modeling accuracy and dispatch optimality of the proposed framework are validated by numerical experiments.

The remaining parts are organized as follows. Section \ref{sec_formulation} formulates the optimal power dispatch model. Section \ref{sec_MLP} presents the construction and reformulation of MLP. Section \ref{sec_case} conducts simulations and Section \ref{sec_conclusion} concludes this paper.

\section{Problem formulation} \label{sec_formulation}

\subsection{Modeling of building thermal dynamics}
A distribution network contains multiple buses and transmission lines.
Each bus can provide electricity for multiple buildings. This paper considers the flexibility of building thermal inertia in the power dispatch. We aggregate all buildings in one bus as a big thermal zone to simplify the problem.\footnote{Note that this will not significantly affect the power dispatch performance and considering multiple buildings under one bus is trivial.} 
We use $i$ and $t$ to index buses/zones ($i \in \mathcal{I} = \{1,2,3,...,I\}$) and time slots ($t \in \mathcal{T} = \{1,2,3,...,T\}$). The thermal dynamics of each aggregated zone can be expressed as follows \cite{hao2016transactive}:
\begin{align}
	\begin{split}
		C_i \frac{d\theta_{i,t}^\text{in}}{dt}=&\frac{\theta_{t}^\text{out}-\theta_{i,t}^\text{in}}{R_i} + q^\text{h}_{i,t} - q^\text{c}_{i,t}, \quad \forall i \in \mathcal{I}, \quad \forall t \in \mathcal{T},
	\end{split} \label{eqn_indoor}
\end{align}
where $\theta_{i,t}^{in}$ and $\theta_{t}^\text{out}$ are the temperatures of indoor environment and ambience, in \textcelsius, respectively. Parameter $R_{i}$ represents the thermal resistances from the $i$-th zone to the surroundings, in \textcelsius/kW. Symbol $q^\text{h}_{i,t}$ and $q^\text{c}_{i,t}$ are heat load and cooling supply of the $i$-th space, in kW, respectively.
Eq. (\ref{eqn_indoor}) can be converted into a linear form by using the finite difference method:
\begin{align}
	\begin{split}
		\theta_{i,t}^\text{in}= &\alpha_{i} \theta_{i,t-1}^\text{in} + \beta_{i} (q^\text{h}_{i,t-1} - q^\text{c}_{i,t-1}) \\
		& + \gamma_{i} \theta_{t-1}^\text{out}, \quad \forall i \in \mathcal{I}, \quad \forall t \in \mathcal{T},
		\label{eqn_T_in}
	\end{split}
\end{align}
where $\alpha_{i}  = 1 - \Delta t/(R_iC_i)$, $\beta_{i} = \Delta t/C_i$ and $\gamma_{i} = \Delta t/(R_iC_i)$. 
The indoor temperature should keep in a proper range to ensure the thermal comforts, as follows:
\begin{align}
	\underline{\theta} \leq \theta_{i,t} \leq \overline{\theta}, \quad \forall i \in \mathcal{I}, \quad \forall t \in \mathcal{T},
\end{align}
where $\underline{\theta}$ and $\overline{\theta}$ denote the lower and upper temperature bounds of thermal comfortable region.

Cooling systems need to consume power to supply cooling power to each zone. The corresponding active power consumption $P_{i,t}^\text{cs}$ in the $i$-th space in time slot $t$ is calculated by:
\begin{align}
	P_{i,t}^\text{cs} = q^\text{c}_{i,t} / COP_i, \quad \forall i \in \mathcal{I}, \quad \forall t \in \mathcal{T},  \label{eqn_P_cs}
\end{align} 
where the parameter $COP$ represents the coefficient of performance of cooling systems. The reactive power demands of cooling systems are set as zeros because modern cooling systems usually contain correction devices to maintain their power factor close to one \cite{sanhueza2007analysis}.

\subsection{Modeling of power flow}
To guarantee the safe operation, all buses and lines have to operate at an appropriate working condition. For example, the bus voltages must keep in [0.9p.u., 1.1p.u.], and the line currents shall not violate allowable upper bounds. 

In order to build the security constraints, we first construct an operation variable vector $\bm x_t \in \mathbb{R}^{3I}$ as:
\begin{align}
	(\bm x_t)^\intercal = [(\bm P_t)^\intercal, (\bm Q_t)^\intercal, (\bm G_t^\text{PV})^\intercal],  \quad \forall t \in \mathcal{T}, \label{eqn_x}
\end{align}
where $\bm P_t$, $\bm Q_t$ and $\bm G_t$ denote the vector form of active, reactive demands and actual used PV generation at all buses in time slot $t$. The active power demands $\bm P_t$ can be expressed as:
\begin{align}
	\bm P_t = \bm P_{t}^\text{cs} + \bm P_{t}^\text{e},\quad \forall t \in \mathcal{T}, \label{eqn_P} 
\end{align}
where $\bm P_{t}^\text{e}$ represents the base load, i.e., the other active power load in the building except that of the cooling system in $t$.\footnote{Note that $\bm Q_t$ is the reactive base load because that we assume the cooling system consume or generate no reactive power.} 

The actually used PV generation should be less than the available PV generation $\bm G_\text{av}^\text{PV}$:
\begin{align}
	\bm G_t^\text{PV} \leq \bm G_\text{av}^\text{PV},\quad \forall t \in \mathcal{T}. \label{eqn_PV} 
\end{align}

When the operation variable vector $\bm x_t$ is determined, the operation state variables, e.g., currents and voltages, can be calculated based on power flow equations. Hence, the security constraints of power flows can be expressed as:
\begin{align}
	\bm s_t(\bm x_t) \in \mathcal{S},\quad \forall t \in \mathcal{T}, \label{eqn_SC}
\end{align}
where $\bm s_t$ is the vector containing all operation state variables at time $t$. Symbol $\mathcal{S}$ denotes the feasible set for safe operation, which is constructed based on security constraints. 

The operation variable vector $\bm x_t$ also determines the power loss. Thus, the power loss $l_t$ is a function of $\bm x_t$:
\begin{align}
 	l_t = f(\bm x_t),\quad \forall t \in \mathcal{T}. \label{eqn_powerloss}
\end{align}

The electricity supply from the main grid $G_t^\text{grid}$ can be calculated via power balance:
\begin{align}
	& G_t^\text{grid} = G_t^\text{buy} - G_t^\text{sell} = \bm 1^\intercal \bm P_t + l_t - \bm 1^\intercal \bm G_t^\text{PV},\quad \forall t \in \mathcal{T},	\label{eqn_balance}\\
	& G_t^\text{buy}\geq 0, G_t^\text{sell}\geq 0, \quad\forall t \in \mathcal{T}. \label{eqn_netload}
\end{align}
If $G_t^\text{grid} = G_t^\text{buy}\geq 0$, then the distribution network purchases power from the main grid; otherwise $G_t^\text{grid} = -G_t^\text{sell}\leq 0$, it sells electricity. The energy cost $EC_t$ can be expressed as:
\begin{align}
	EC_t = c^\text{buy} G_t^\text{buy} - c^\text{sell} G_t^\text{sell},\quad \forall t \in \mathcal{T}, \label{eqn_EC}
\end{align}
where {$c^\text{buy}$ and $c^\text{sell}$} denotes the prices for purchasing and selling electricity. 
\subsection{Formulation of optimal power dispatch problem}
We propose a security-constrained optimal power dispatch strategy for distribution networks. The building thermal inertia is treated as energy storage systems to provide flexibility. The proposed model can be summarized as follows:
\begin{align} 
&\min_{\{q_t^\text{t}, \forall i \in \mathcal{I}, \forall t \in \mathcal{T}\}} \quad \sum_{t=1}^{T} EC_t \tag{$\mathcal{P}_1$}, \quad
\text{s.t.:}\quad \text{Eqs. (\ref{eqn_T_in})-(\ref{eqn_EC})}.
\notag
\end{align}
Note that the decision variables in the dispatch model include cooling supplies $q_t^\text{t}$ and actual PV generation $\bm G_t^\text{PV}$; The other parts in the operation vector $\bm x_t$ are treated as given parameters.

\begin{remark}
	Conventional methods require accurate power flow models to describe function $\bm s_t(\bm x_t)$ in Eq. (\ref{eqn_SC}) and $f(\bm x_t)$ in Eq. (\ref{eqn_powerloss}). However, in many practical cases, an accurate power flow model may be unavailable because of unclear topology information or unknown network parameters. 
\end{remark}

\section{Topology-free optimal power dispatch} \label{sec_MLP}
To realize the topology-free power dispatch framework, a multi-layer perception (MLP) neural network is trained to classify whether the state variables $\bm s_t$ stay in the feasible set $\mathcal{S}$ with any given operation vector $\bm x_t$. 
Then, the forward propagation of MLP is exactly reformulated into  mixed-integer constraints that are used to substitute constraint (\ref{eqn_SC}) in problem $\mathcal{P}_1$. A linear regression (LR) is also trained to approximate the power loss function (\ref{eqn_powerloss}).  As a result, the optimal dispatch strategy becomes a mixed integer linear program, which can be efficiently solved by off-the-shelf solvers. 
\subsection{Construction of MLP}
\begin{figure}
	\centering
				\vspace{-4mm}
	\includegraphics[width=0.8\columnwidth]{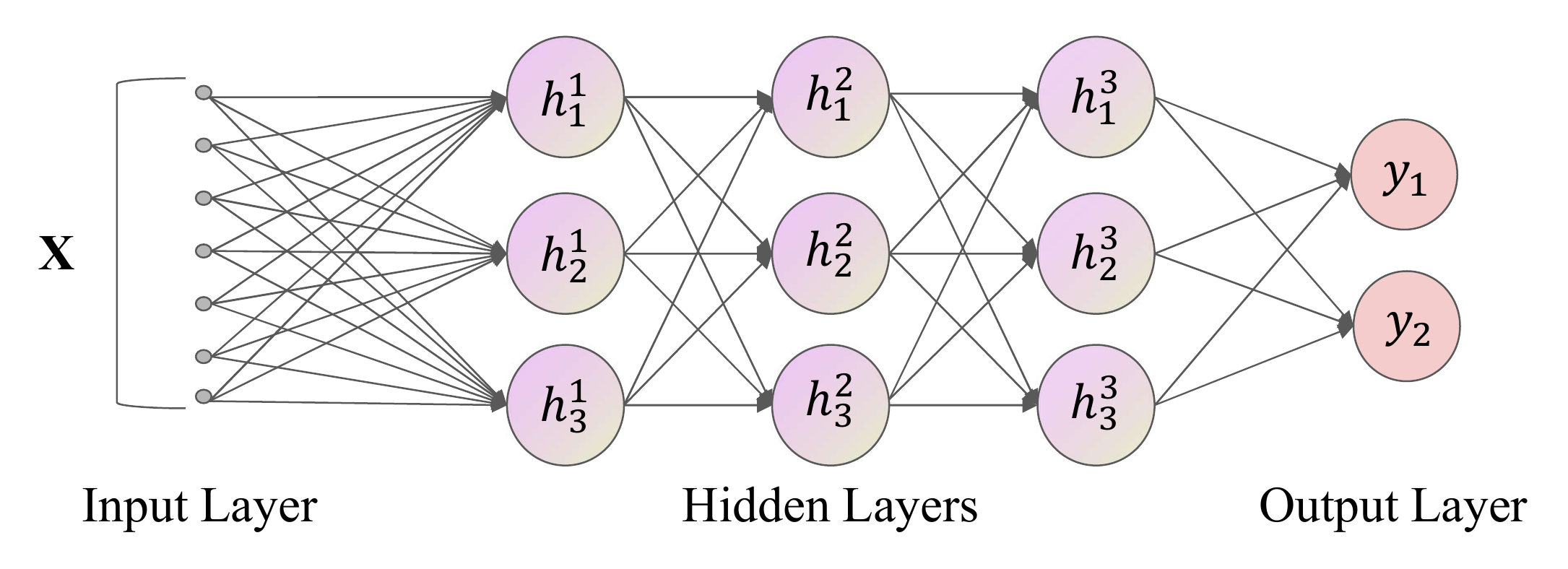}\vspace{-2mm}
	\caption{A MLP with 3 hidden layers, where $h_n^k$ denotes the output of $n$-th neuron in the $k$-th hidden layer.
	}
	\label{fig_MLP}
	\vspace{-4mm}
\end{figure}
A MLP is trained to judge the safety of current operation variable $\bm x_t$. Fig. \ref{fig_MLP} shows the typical structure of MLP. 
It consists of one input layer, $K$ hidden layers and one output layer. Each layer contains multiple neurons. 
Every neuron is composed of a linear mapping and a non-linear activation function. We use $k$ and $n_k$ to index the hidden layers ($k \in \mathcal{K} = {1,\cdots, K}$) and the neurons in the $k$-th hidden layer ($n_k \in \mathcal{N}_k = {1,\cdots, N_k}$). ReLU is employed as the non-linear activation function. Then, the forward propagation of MLP is expressed as
\begin{align}
    &\bm x_t = \bm h^{0}_t, \quad \forall t \in \mathcal{T}, \label{eqn_input}\\
	&\bm z^{k}_t = \bm W^{k} \bm h^{k-1}_t + \bm b^{k}, \quad \forall k \in \mathcal{K}, \quad \forall t \in \mathcal{T}, \label{eqn_z} \\
	&\bm h^{k}_t = \max(\bm z^{k}_t, 0), \quad \forall k \in \mathcal{K},\quad \forall t \in \mathcal{T}, \label{eqn_h} 
\end{align}
where $h^k_t \in \mathbb{R}^{N_k}$ represents the outputs of layer $k$. Parameters $\bm W^k \in \mathbb{R}^{N_{k} \times N_{k-1}}$ and $\bm b^{k} \in \mathbb{R}^{N_{k}}$ denote the weight matrix and bias of layer $k$, which are parameters to be trained. 
The outputs of this MLP, i.e., $\bm y_t \in \mathbb{R}^{2}$, is given by:
\begin{align}
	\bm y_t = \bm W^{K+1} \bm h^{K}_t + \bm b^{K+1},\quad \forall t \in \mathcal{T}. \label{eqn_output}
\end{align}
The final classification result is determined by comparing the two elements of $\bm y_t$. The sample is labeled as ``unsafe" if $y_{1,t}>y_{2,t}$; otherwise, this sample is labeled as ``safe".

\subsection{Construction of LR model}
In distribution networks, the effects of power loss on the total energy cost is not significant because it is usually small compared with the total active demands. It is not necessary to train a time-consuming model for predicting power loss with high accuracy. Thus, a simple LR is used to build the mapping between operation variable $\bm x_t$ and power loss $l_t$, as follows:
\begin{align}
l_t = \bm W^\text{loss} \bm x_t + b^\text{loss},\quad \forall t \in \mathcal{T}, \label{eqn_lr}
\end{align}
where $W^\text{loss} \in \mathbb{R}^{3I}$ and $b^\text{loss} \in \mathbb{R}$ are parameters to be trained.

 
\subsection{Reformation of topology-free power dispatch model}
The constraints (\ref{eqn_z})-(\ref{eqn_h}) are intractable because of the maximum operator. 
However, they can be reformulated into mixed-integer forms by introducing auxiliary continuous variables $\bm r^{k}_t$ and binary variable $\bm \mu^{k}_t$ based on the big-M method \cite{fischetti2018deep}:
\begin{align}
&\bm h^{k}_t - \bm r^{k}_t = \bm W^{k} \bm h^{k-1}_t + \bm b^{k}, \quad \forall k \in \mathcal{K}, \quad \forall t \in \mathcal{T}, \label{eqn_r1} \\
&0 \leq \bm h^{k}_t \leq M \cdot \bm \mu^{k}_t, \quad \forall k \in \mathcal{K},  \quad \forall t \in \mathcal{T}, \label{eqn_r2} \\
&0 \leq \bm r^{k}_t \leq M\cdot(1-\bm \mu^{k}_t), \quad \forall k \in \mathcal{K}, \quad \forall t \in \mathcal{T}, \label{eqn_r3} \\
&\bm \mu^k_t \in \{0,1\}^{N_k}, \quad \forall k \in \mathcal{K},  \quad \forall t \in \mathcal{T}, \label{eqn_r4}
\end{align}
in which $M$ is a big enough constant. 

By replacing the security constraint (\ref{eqn_SC}) and the power loss function (\ref{eqn_powerloss}) with the reformulated MLP (\ref{eqn_r1})-(\ref{eqn_r4}) and LR model (\ref{eqn_lr}), respectively, problem $\mathcal{P}_1$ can be reformulated as:
\begin{align} 
	&\min_{\{q_t^\text{t}, \forall i \in \mathcal{I}, \forall t \in \mathcal{T}\}} \quad \sum_{t=1}^{T} EC_t \tag{$\mathcal{P}_2$}\\
	&\begin{array}{r@{\quad}r@{}l@{\quad}l}
		\text{s.t.:} & &\text{Eqs. (\ref{eqn_T_in})-(\ref{eqn_P}), (\ref{eqn_balance})-(\ref{eqn_input}), (\ref{eqn_lr})-(\ref{eqn_r4})}, \\
		&& y_{t,1} \leq y_{t,2}, \quad \forall t \in \mathcal{T}.
	\end{array} \notag
\end{align}
In $\mathcal{P}_2$, constraints (\ref{eqn_T_in})-(\ref{eqn_P_cs}) describe the thermal dynamics of buildings. Constraints (\ref{eqn_x})-(\ref{eqn_P}) give the relation between original decision variables and inputs of MLP and LR.  Constraints (\ref{eqn_balance})-(\ref{eqn_EC}) define the power balance and calculate the energy cost in objective function. Constraint (\ref{eqn_input}) defines the inputs of MLP. Constraint (\ref{eqn_lr}) calculates the power loss. Constraints (\ref{eqn_r1})-(\ref{eqn_r4}) represent the forward propagation of MLP. The last constraint $y_{t,1} \leq y_{t,2}$ ensures that the obtained solution satisfies the security constraints of power flow. 

\begin{remark}
	The proposed framework only needs historical data to train the MLP and LR. Thus, the security constraints and power loss can be introduced to our optimization without any topology information.
	Besides, $\mathcal{P}_2$ is a mixed-integer linear program. It can be efficiently solved by the off-the-shelf solvers. The global optimality can be also guaranteed.  
\end{remark}

\section{Case study} \label{sec_case}
\begin{figure}
	\centering
	\vspace{-8mm}
	\includegraphics[width=0.8\columnwidth]{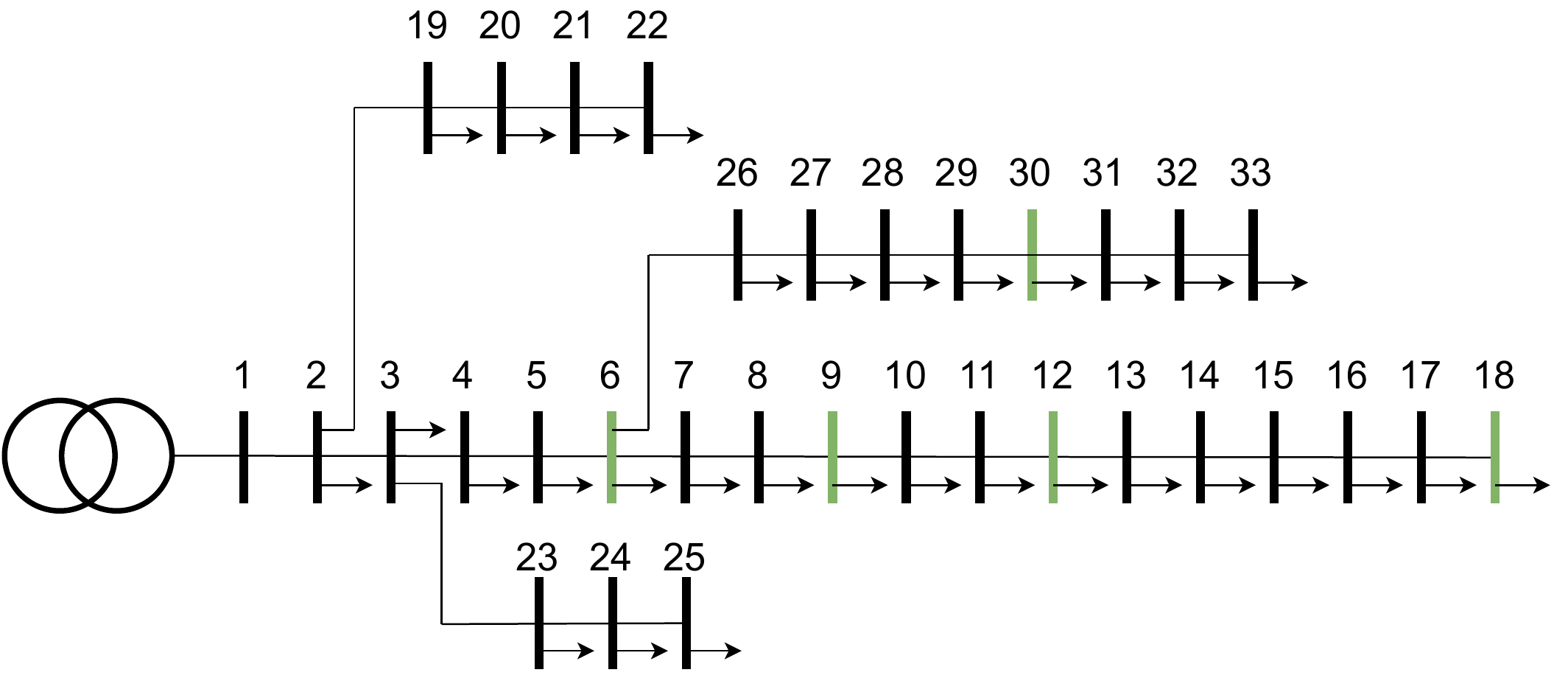}\vspace{-2mm}
	\caption{Structure of the 33-bus system. PV power stations are installed at green buses.
	}
	\label{fig_33bus}
	\vspace{-4mm}
\end{figure}

\begin{figure}
	\centering
	\subfigure[]{\includegraphics[width=0.49\columnwidth]{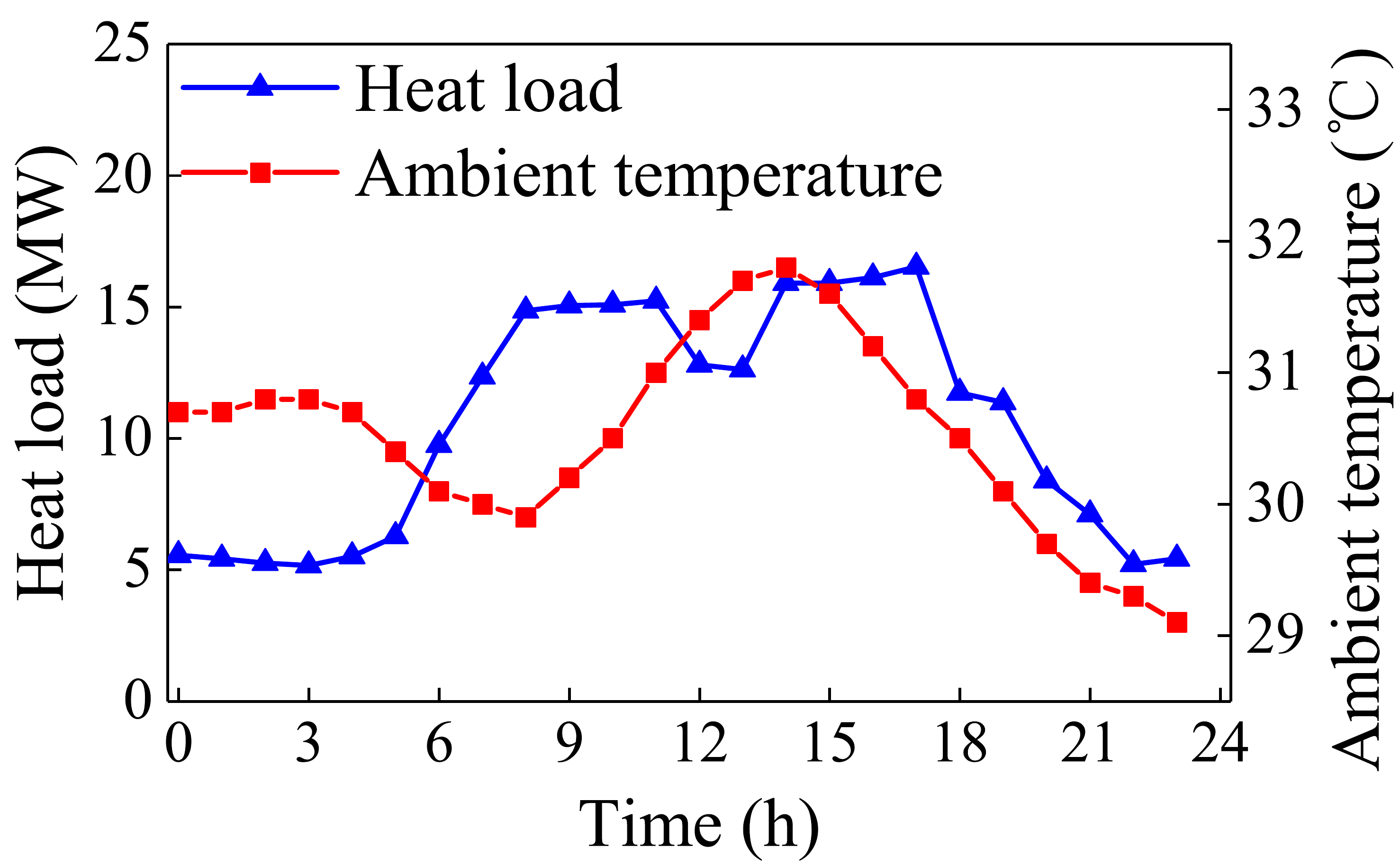}}
	\subfigure[]{\includegraphics[width=0.49\columnwidth]{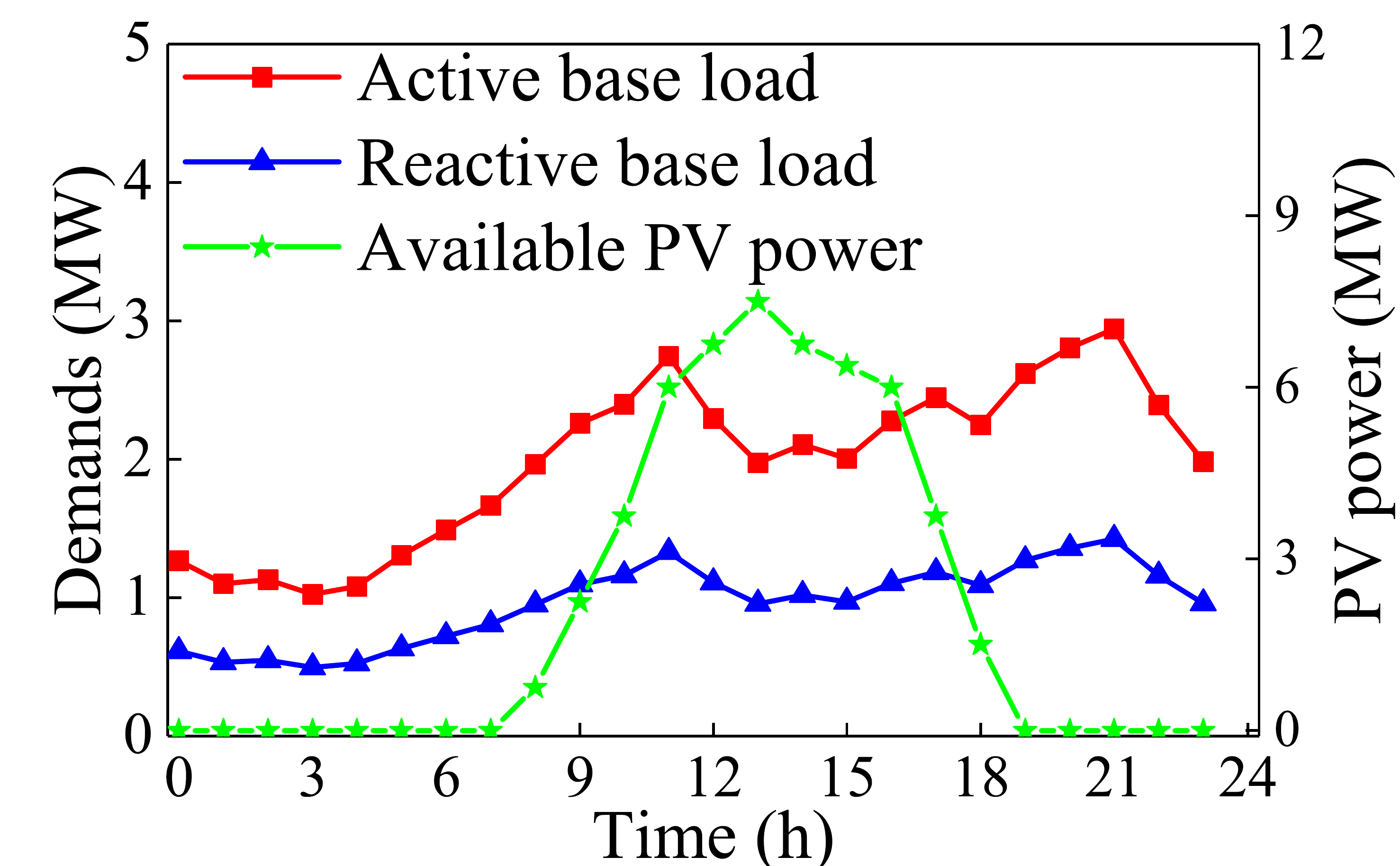}}\vspace{-3mm}
	\vspace{-2mm}
	\caption{(a) Total heat loads, ambient temperature, and (b) active based load, reactive base load and available PV generation.}
	\label{fig_parameters}
	\vspace{-4mm}
\end{figure}

We implement two case studies (i.e. heavy and light load cases) based on the IEEE 33-bus system \cite{baran1989network} (shown in Fig. \ref{fig_33bus}) to demonstrate the benefits of the proposed framework. The heavy-load case is used to validate the optimality and reliability of the proposed MLP-based framework, while the light-load case is implemented to show the benefits of treating building thermal inertia as flexibility source. We assume PV power stations are installed at 5 different buses (i.e. buses 6, 9, 12, 18 and 30). Note that the information of topology and parameters of the IEEE 33-bus system are only used to generate training data and assumed to be unknown in the view of the proposed framework. {The total heat load $\bm q_t^h$, ambient temperature $\theta_t^\text{out}$, active based load $\bm P_t^e$, reactive base load $\bm Q_t$ and available PV generation $\bm G_\text{av}^\text{PV}$ are illustrated in Fig. \ref{fig_parameters}.} All loads in the heavy-load case are the same with those in Fig.\ref{fig_parameters}, while the loads are reduced by half in the light-load case.
The optimization time horizon is 24 hours. The values of other parameters are summarized in Table \ref{tab_parameter}.

\begin{table}
	\small
	\centering
		\vspace{-4mm}
	\caption{Parameters in case study}
		\vspace{-2mm}
	\begin{tabular}{cccc}
		\hline
		\rule{0pt}{11pt}		
		Parameters & Value &Parameters & Value\\
		\hline
		\rule{0pt}{10pt}
		$C_i$ &  1 MWh/{\textcelsius} & $\underline \theta$& 24{\textcelsius}\\
		$R_i$ &  50 {\textcelsius}/MW & $\overline \theta$& 28{\textcelsius} \\
		$COP_i$ &  3.6 & $\eta^\text{buy}$ & 0.1122 \$/kWh \\
		$\Delta t$ & 1h & $\eta^\text{sell}$ & 0.056 \$/kWh \\		
		\hline
	\end{tabular}\label{tab_parameter}
	\vspace{-4mm}
\end{table}

To create the labeled dataset for MLP and LR, we randomly generate the operation vector $\bm x_t$ as features and input them to Pandapower, a power system modeling tool in Python, to calculate state variables, i.g., nodal voltages and distribution line currents, to obtain the corresponding labels. If the vector $\bm x_t$ satisfies all security constraints, this sample's label is set as ``safe"; otherwise, its label is ``unsafe". 
In this paper, the voltage limitation (i.e. the bus voltages are limited in [0.9 p.u., 1.1 p.u.]) and current limitation (i.e. the line currents should be less than 0.249kA) are considered as our security constraints. Finally, 10,000 samples are generated, in which 40\% of them are ``safe"; 60\% of them are ``unsafe". Then, 70\% of samples are randomly selected to form the training set and the remainings are used as the testing set. 

A model without considering the security constraints is implemented as benchmark 1 to show the reliability of the proposed $\mathcal{P}_2$. Benchmark 1 has the same objective function with $\mathcal{P}_2$ but its constraints only contain Eqs. (\ref{eqn_T_in}-\ref{eqn_PV}), (\ref{eqn_balance}-\ref{eqn_EC}) and (\ref{eqn_lr}). 
%
To demonstrate the optimality of the proposed framework, a conventional power-flow-based model is also presented as benchmark 2. This benchmark is a second-order cone program and requires exact topology parameters to build the power flow equations in order to describe security constraints and calculate power losses.

All simulations are performed on an Intel(R) Core(TM) 8700 3.20GHz CPU with 16 GB memory. The MLP and LR are implemented by Pytorch. The corresponding optimization problem is built by CVXPY and solved by GUROBI. 

\subsection{Heavy load case}

\subsubsection{Reliability}
\begin{figure}
	\centering
		\vspace{-4mm}
	\subfigure[]{\includegraphics[width=0.49\columnwidth]{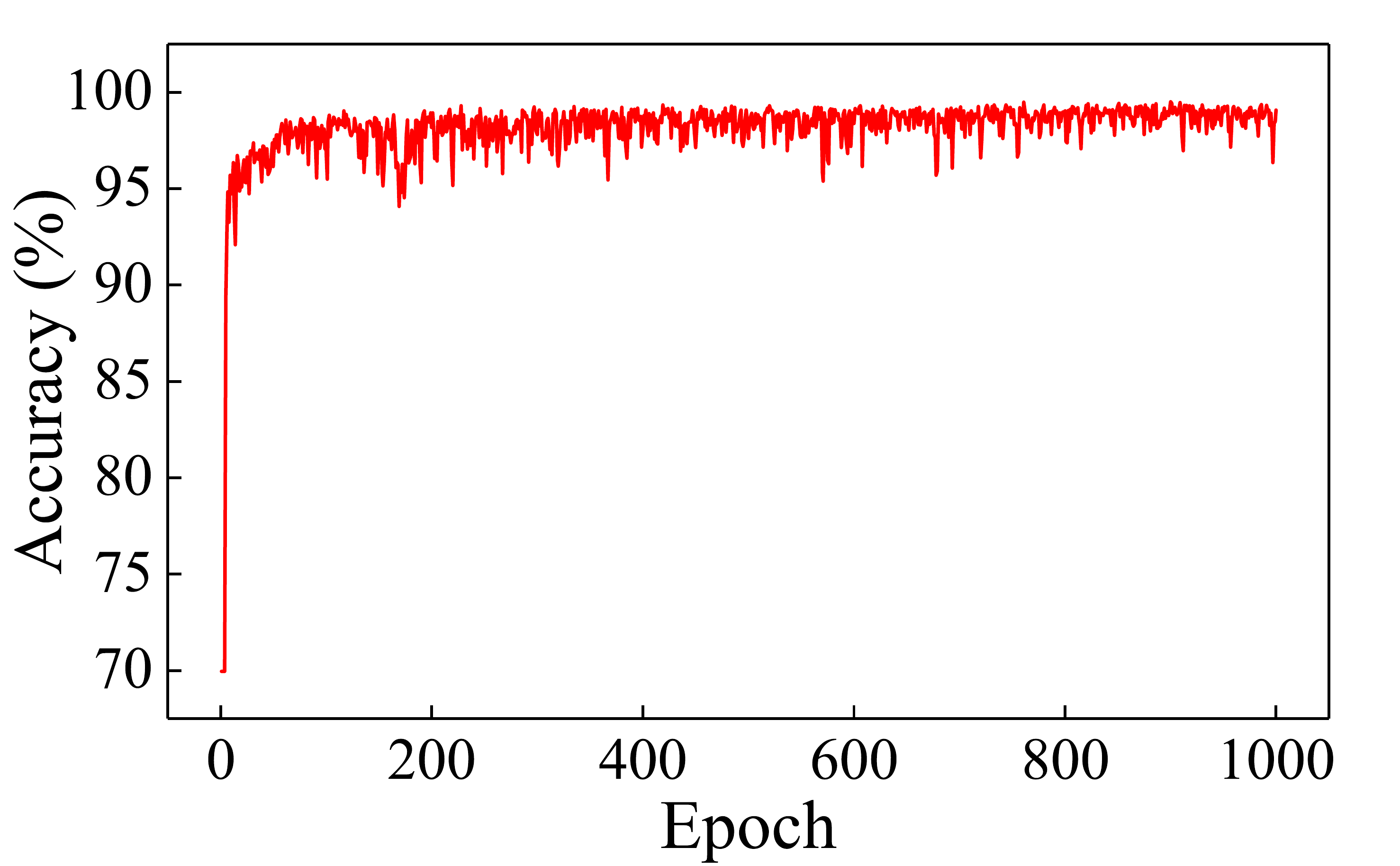}}
	\subfigure[]{\includegraphics[width=0.49\columnwidth]{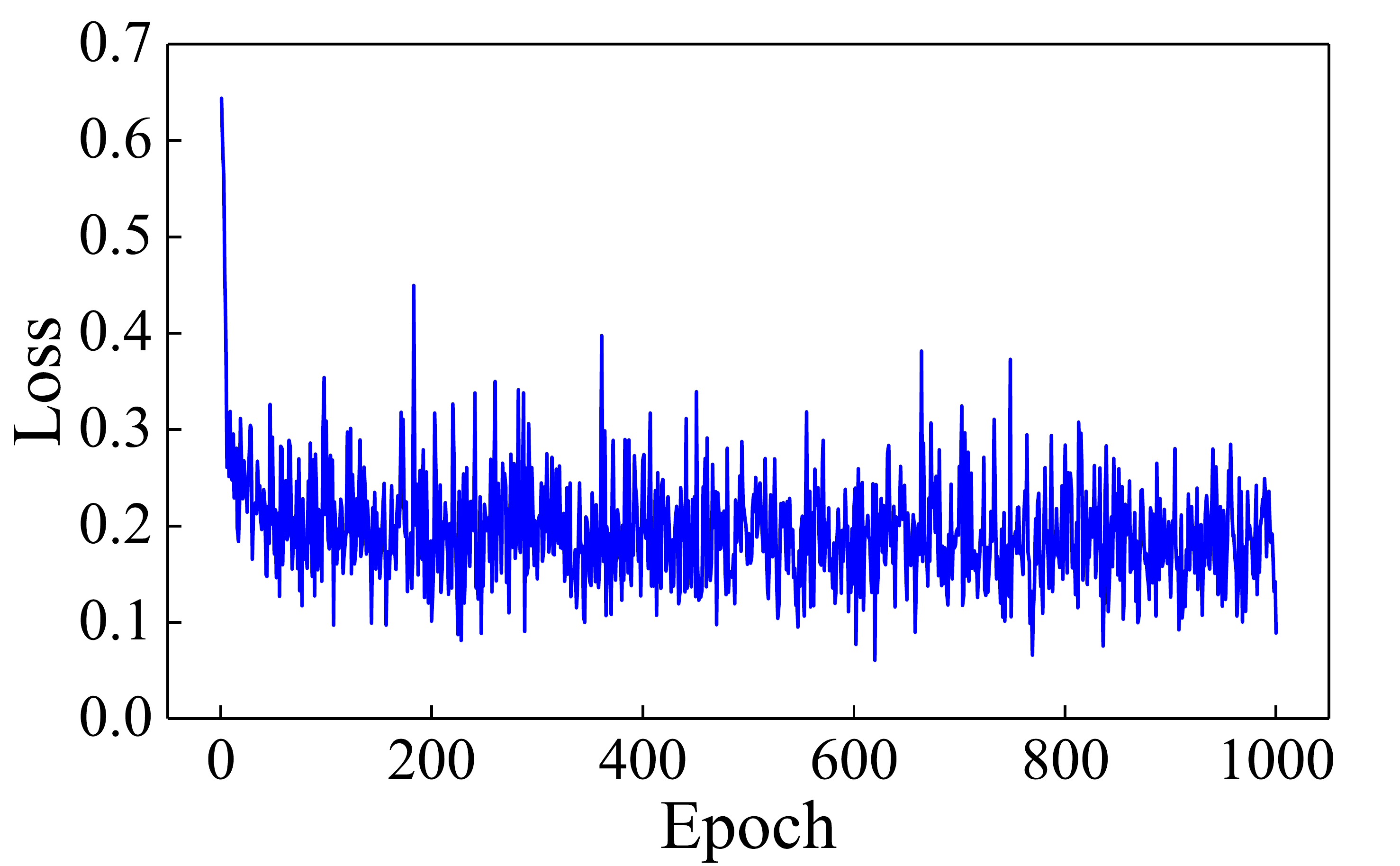}}\vspace{-3mm}
	\vspace{-2mm}
	\caption{(a) Testing accuracy, and (b) {value of loss function during training}.}
	\label{fig_training}
	\vspace{-4mm}
\end{figure}

We use a MLP with 2 hidden layer with 8 neurons in each hidden layer to judge whether the current operation variables $\bm x_t$ stratifies the security constraint (\ref{eqn_SC}). The accuracy on test set and loss function value during training are shown in Fig. \ref{fig_training}. The model accuracy in the testing set can reach to 99.09\% with a low loss function value, which indicates the high accuracy of the MLP. 

\begin{figure}
	\centering
	\subfigure[]{\includegraphics[width=0.49\columnwidth]{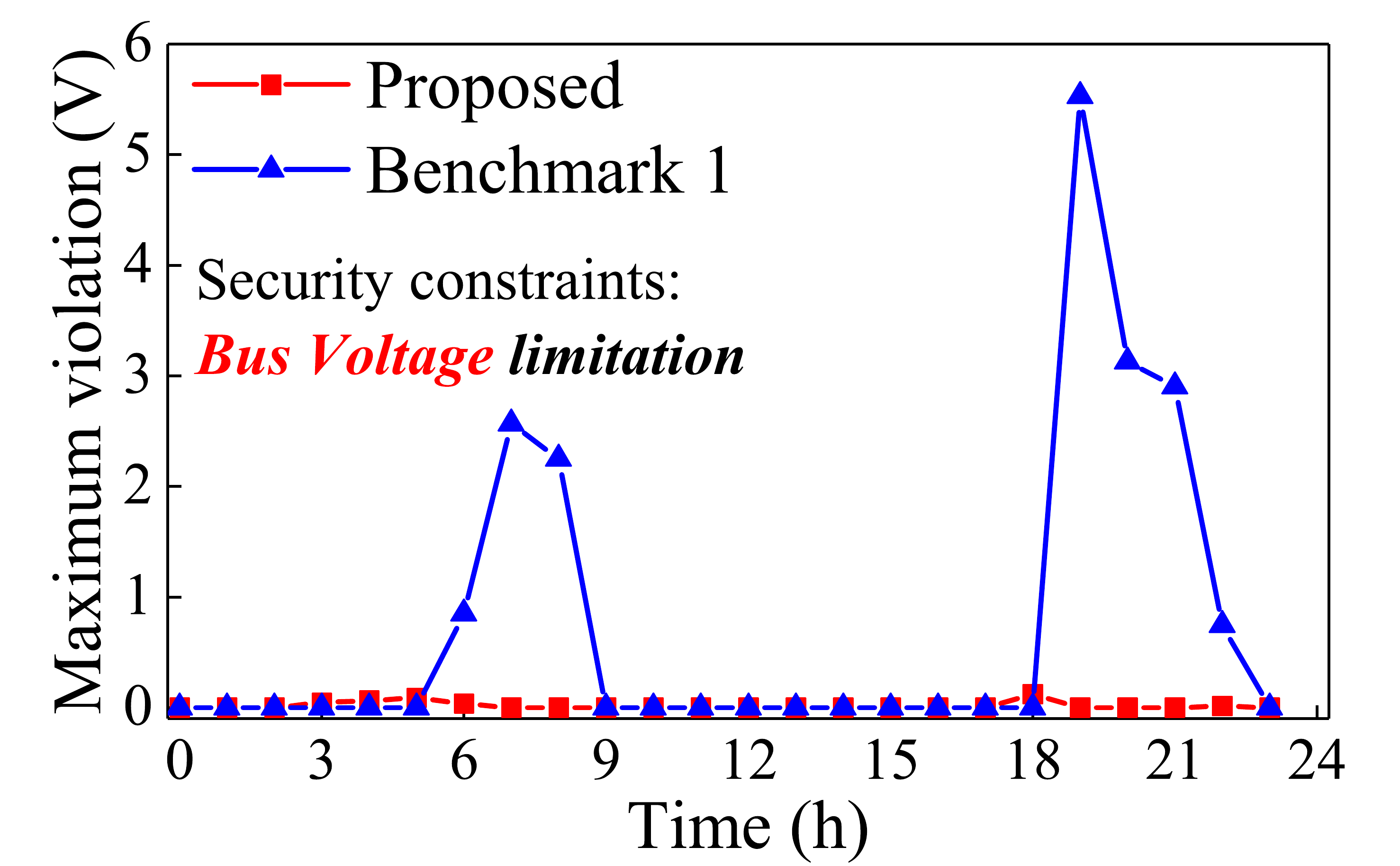}}
	\subfigure[]{\includegraphics[width=0.49\columnwidth]{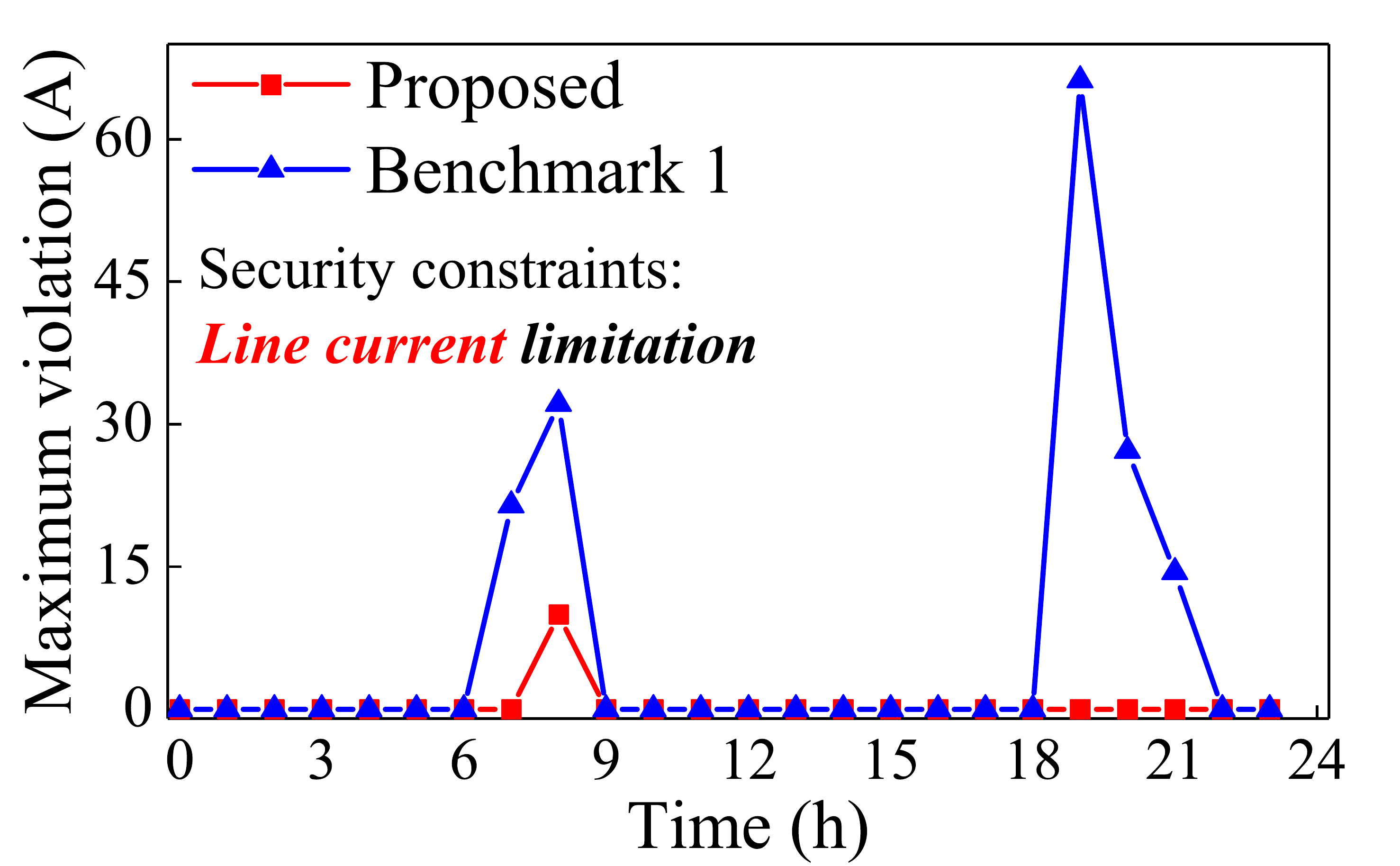}}\vspace{-3mm}
	\vspace{-2mm}
	\caption{Maximum violations of (a) bus voltages, and (b) line currents.}
	\label{fig_violation}
	\vspace{-4mm}
\end{figure}

To demonstrate the reliability of the proposed framework, 
we calculate the corresponding true currents and voltages with Pandapower by treating the solutions of $\mathcal{P}_2$ and benchmark 1 as inputs. The corresponding maximum violations of security constraints among all buses and lines are shown in Fig. \ref{fig_violation}.
The proposed framework can obtain relatively safe strategy, while the solution derived by {benchmark 1} is much worse. For example, the maximum voltage violation of the proposed framework almost keeps at zero. The current violation only appears at $t=8$ and the violation value is 9.99A. The maximum violation of benchmark 1 are 5.53V ($t=19$) and 66.19A ($t=19$), which are much larger than those in the proposed framework. The violations in benchmark 1 also occurs much more frequently. These results confirm the proposed framework can drive a high reliable strategy without any topology information.

\subsubsection{Optimality}
\begin{figure}
	\centering
				\vspace{-4mm}
	\includegraphics[width=0.8\columnwidth]{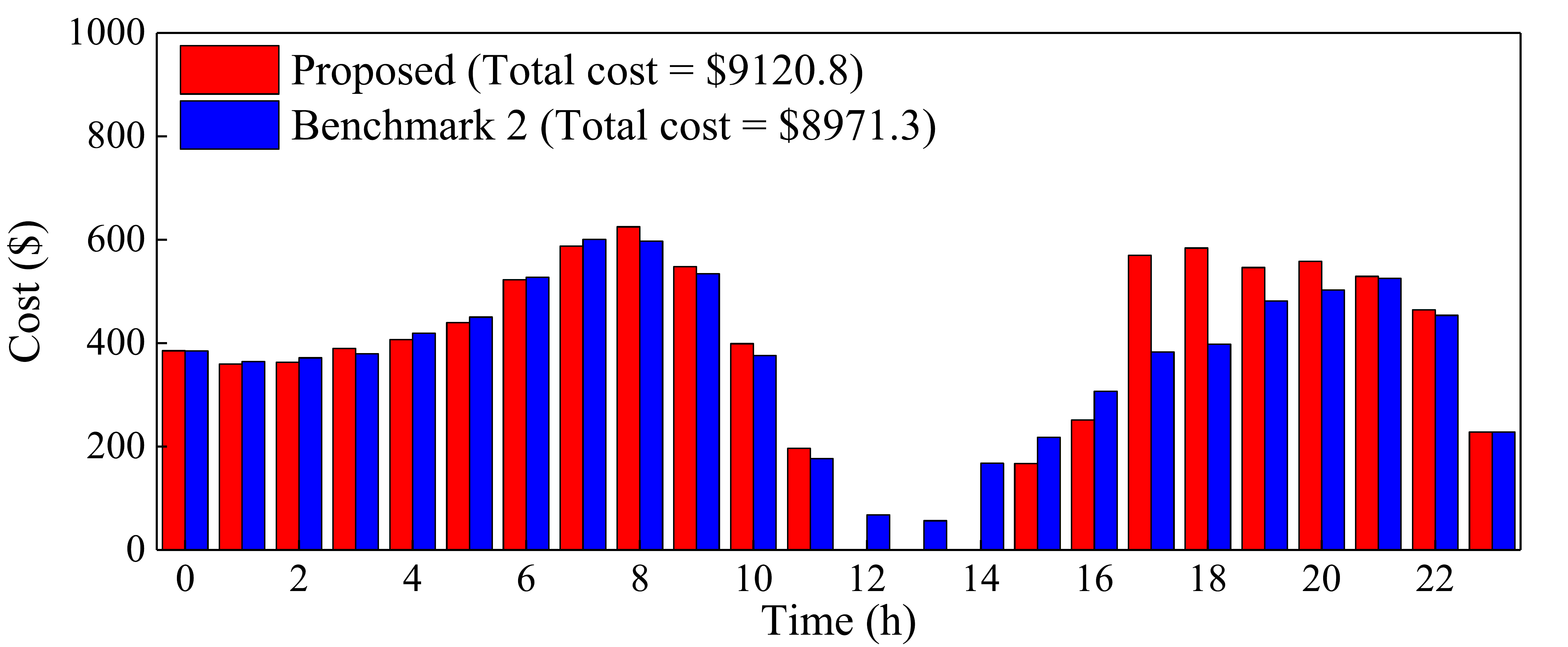}\vspace{-2mm}
	\caption{Comparison between the hourly costs obtained by the proposed framework and benchmark 2.
	}
	\label{fig_hourlyCost}
	\vspace{-4mm}
\end{figure}

Fig. \ref{fig_hourlyCost} compares the hourly costs obtained by $\mathcal{P}_2$ and benchmark 2 to demonstrate the optimality of the proposed framework. In most hours, the hourly cost obtained by the proposed framework is almost the same with that in benchmark 2. During 12:00-16:00, the pre-cooling in the proposed model is not so sufficient like that in benchmark 2, so it has to pay more in the next few hours. Nevertheless, the total cost of the proposed framework is \$9120.8, which is only 1.6\% larger than that in the benchmark 2. These results indicate the high optimality of the proposed framework.

\subsection{Light load case}
We employ the proposed framework to solve the light-load case to demonstrate the benefits of employing building thermal inertia to provide flexibility.
\subsubsection{Indoor temperature}
\begin{figure}
	\centering
				\vspace{-4mm}
	\includegraphics[width=0.8\columnwidth]{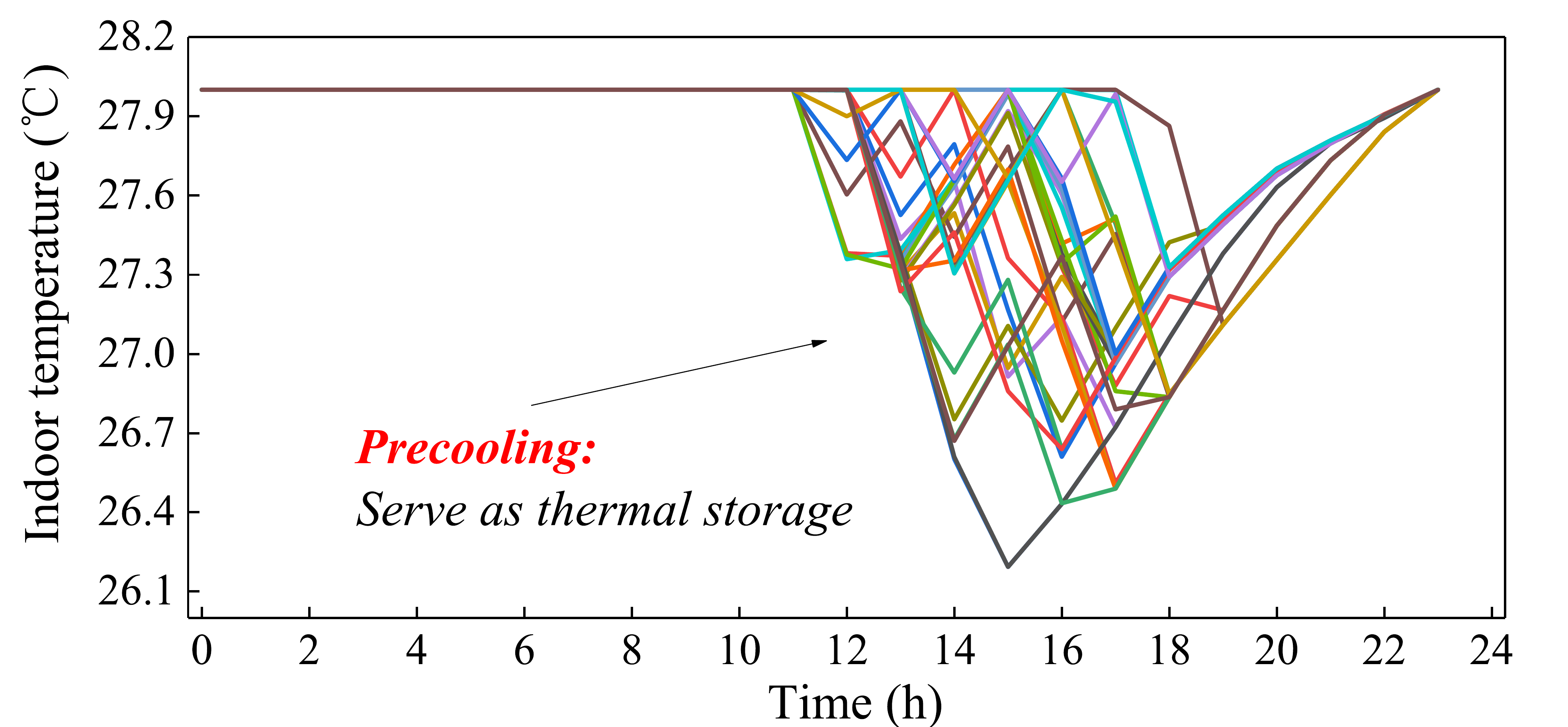}\vspace{-2mm}
	\caption{Results of indoor temperatures with building thermal inertia.
	}
	\label{fig_T_in}
	\vspace{-4mm}
\end{figure}

Fig. \ref{fig_T_in} illustrates the average indoor temperatures obtained by the proposed $\mathcal{P}_2$. In the hours with low PV generation (i.e. 00:00am-11:00am), all indoor temperature keep at the upper bound of the comfortable region to minimize the heating effects from outdoor to indoor environments. The temperature drop occurs at the time with high PV generation (i.e. 11:00am-18:00pm). We call this temperature drop as pre-cooling.  With the decrease of the indoor temperature, more and more power can be stored in the indoor environments by pre-cooling to compensate the cooling demands in the next few hours.

\subsubsection{PV curtailments}
\begin{figure}
	\centering
		\vspace{-2mm}
	\subfigure[]{\includegraphics[width=0.49\columnwidth]{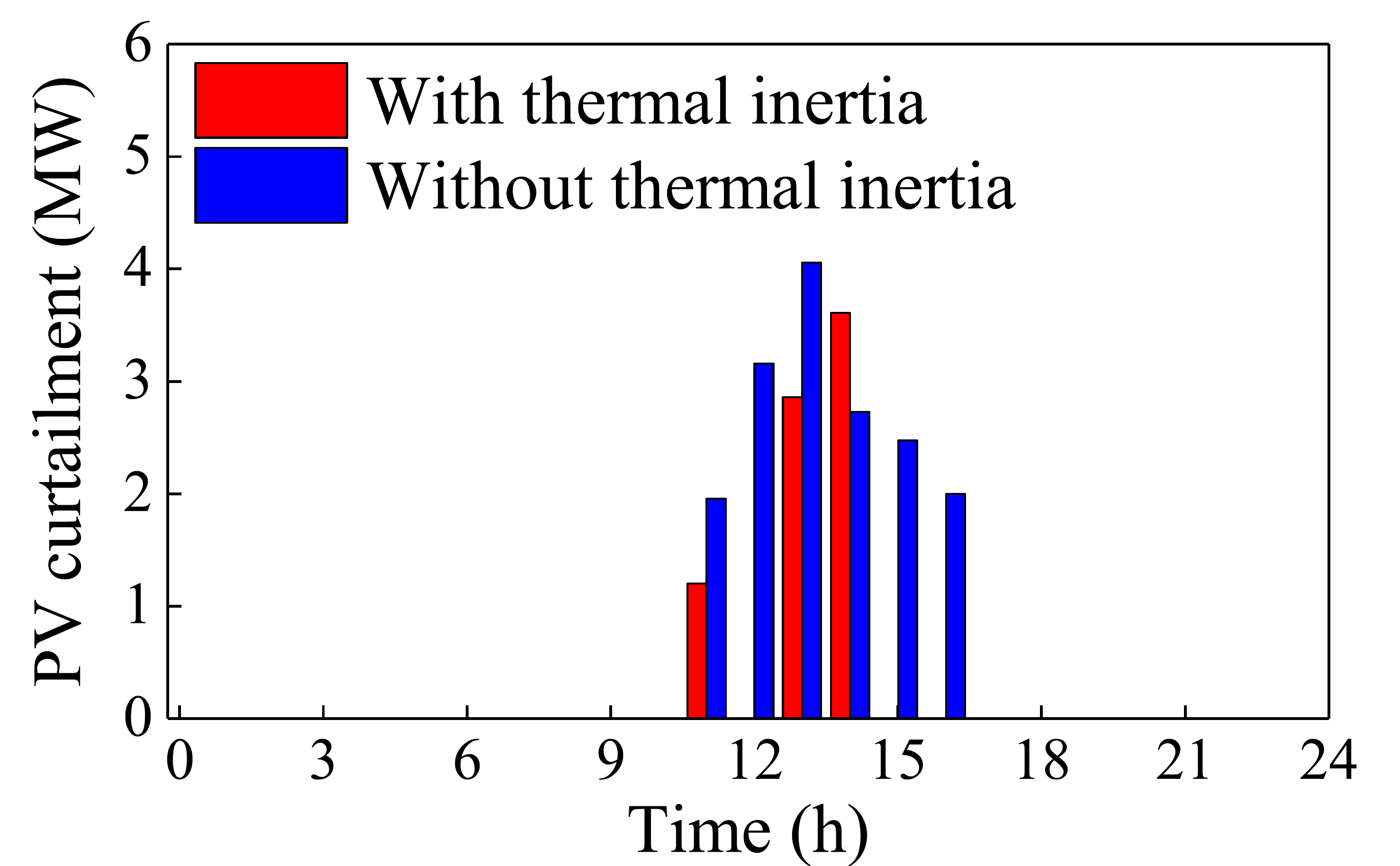}}
	\subfigure[]{\includegraphics[width=0.49\columnwidth]{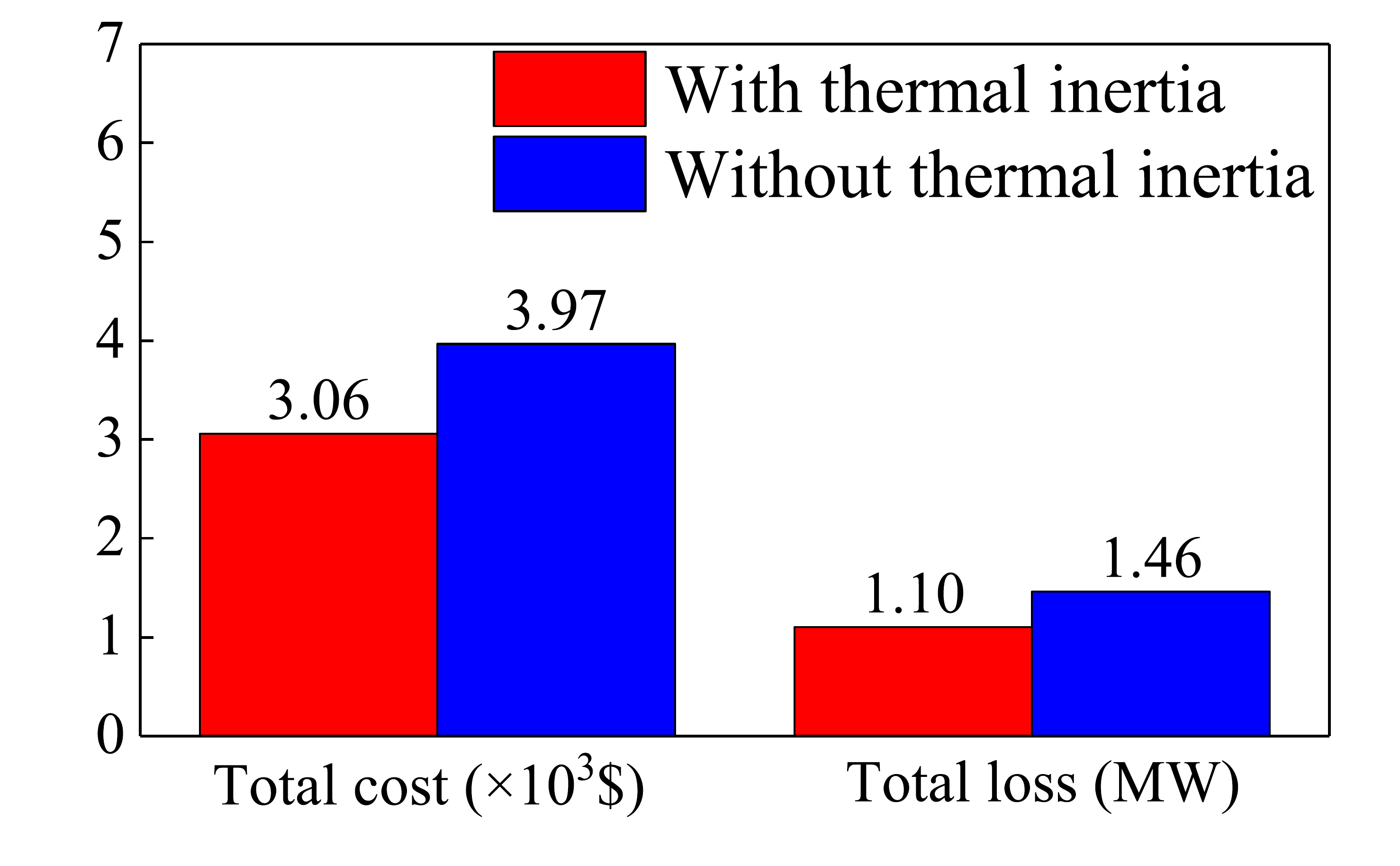}}\vspace{-3mm}
	\vspace{-2mm}
	\caption{Comparison of (a) PV curtailments and (b) total cost and loss with and without the utilization of building thermal inertia.}
	\label{fig_PV_curtailment}
	\vspace{-4mm}
\end{figure}

Fig. \ref{fig_PV_curtailment} compares the results with and without the utilization of building thermal inertia. The scenario without thermal inertia is implemented by restricting the indoor temperature at its maximum value $\overline \theta$. Because the indoor environments can serve as thermal storage system, the excess PV generation can be stored by pre-cooling at the time with high available PV power. Then, the PV curtailment can be significantly reduced, as shown in Fig. \ref{fig_PV_curtailment}(a). 
Meanwhile, because of the thermal inertia, the stored cooling power can satisfy parts of cooling demands in the next few hours without sunshine (i.e. 18:00pm-24:00pm). Thus, the electricity consumed by cooling systems in those hours will also decrease. Besides, lower power consumption means smaller power flow, which can also reduce power loss.

\section{Conclusions} \label{sec_conclusion}
This paper proposes a topology-free framework to realize the security-constrained optimal power dispatch for distribution networks. The building thermal inertia is employed as thermal storage systems to store the extra PV generation. Unlike the conventional power-flow-based approach which requires accurate topology information, the proposed framework trains a MLP to replace the security constraints. To achieve the global optima efficiently, the MLP is reformulated as new constraints with mixed-integer form. Simulation results show that the classification accuracy of MLP can achieve 99.09\%. Then, two case studies are conducted to demonstrate the feasibility, optimality of the proposed framework and the benefits of utilizing building thermal inertia. The results show that the proposed framework can obtain a proper operation strategy with fewer violations of security constraints. Moreover, the relative optimality gap is only 1.6\%. The results also indicate that employing thermal inertia can effectively reduce PV curtailment and power losses, and save total energy cost.

\footnotesize
\bibliographystyle{ieeetr}
\bibliography{ref}

\appendices
\setcounter{table}{0}   
%

\end{document}